\documentclass[twoside, 12pt]{amsart}
\usepackage{latexsym}
\usepackage{amssymb,amsmath,amsopn}
\usepackage[dvips]{graphicx}   
\usepackage{color,epsfig}      

\newtheorem{theorem}{Theorem}
\newtheorem{lemma}{Lemma}
\newtheorem{proposition}{Proposition}

\newtheorem{remark}{Remark}

\newtheorem{conjecture}{Conjecture}
\newtheorem{question}{Question}

\newcommand{\ee}{\mathbf{e}}
\newcommand{\vv}{\mathbf{v}}
\newcommand{\uu}{\mathbf{u}}

\DeclareMathOperator{\grad}{grad}
\DeclareMathOperator{\dif}{d}
\DeclareMathOperator{\vol}{vol}
\DeclareMathOperator{\surf}{surf}
\DeclareMathOperator{\area}{area}

\DeclareMathOperator{\bd}{bd}

\renewcommand{\Re}{\mathbb R}
\newcommand{\Sph}{\mathbb{S}}

\newcommand{\B}{\mathbf B}

\begin{document}

\title[Shape descriptors evolving under distance-driven flow]{The evolution of geophysical shape descriptors under distance-driven flows}
\author[G. Domokos and Z. L\'angi]{G\'abor Domokos and Zsolt L\'angi}
\address{G\'abor Domokos, MTA-BME Morphodynamics Research Group and Dept. of Mechanics, Materials and Structures, Budapest University of Technology,
M\H uegyetem rakpart 1-3., Budapest, Hungary, 1111}
\email{domokos@iit.bme.hu}
\address{Zsolt L\'angi, MTA-BME Morphodynamics Research Group and Dept. of Geometry, Budapest University of Technology,
Egry J\'ozsef u. 1., Budapest, Hungary, 1111}
\email{zlangi@math.bme.hu}

\begin{abstract}
We investigate the evolution of axis ratios, roundness (isoperimetric ratio)  and the number of static balance points under distance-driven flows.
The latter have already been proposed by Aristotle as models of particle shape evolution and recent studies indicate that they may serve as models for frictional abrasion.
We show exact conditions under which Aristotle's original claims are true.  For several  geophysical shape descriptors we prove monotonic or quasiconcave time evolution and
compare these results with results from the literature on curvature-driven flows as models of collisional abrasion.
\end{abstract}

\maketitle


\section{Introduction}

Physical abrasion processes, composed of collisional and frictional abrasion, are of fundamental importance in the evolution of sedimentary particles.
Geologists try to track the shape evolution process by measuring scalar quantities called shape descriptors associated with the particle's shape.
The most common shape descriptors are the axis ratios of the approximating ellipsoid (Zingg 1935) and roundness (Cox 1927) of the orthogonal projection (sometimes
also referred to as circularity (Blott and Pye 2008)). More recently, the number of static balance points has been introduced (Domokos et al. 2010) as a useful shape descriptor.
While substantial amount of data on shape descriptors has accumulated over decades (e.g. Bluck 1967; Carr 1969, 1972; Griffith 1967; Zingg 1935), and there is growing literature on curvature-driven flows (Bloore 1977; Domokos 2014; Domokos and Gibbons 2012, 2013; Firey 1974; Miller et al. 2014) serving as mathematical models of collisional abrasion, very little is known on the evolution of shape descriptors under frictional abrasion.  While the mathematical models of the latter still lack rigorous experimental verification, the framework proposed in (Domokos and Gibbons 2012, 2013) for frictional abrasion suggests that distance-driven flows may be the best candidate models.
These flows also have great historic importance as the Aristotelian models of particle abrasion (Domokos and Gibbons 2012; Krynine 1960).

In the current paper we investigate the evolution of shape descriptors under distance-driven flows with special emphasis on Aristotelian models (radial flows) and potential models of frictional abrasion (parallel flows).
The main goal of our work is to find cases when the time evolution of a given shape descriptor is simple from the geophysical point of view.
In mathematical terms this means that the shape descriptor evolves either monotonically (its time evolution has no extremum), in a quasiconcave manner (its time evolution
has one single maximum and no minima) or in a quasiconvex manner (its time evolution
has one single minimum and no maxima). Previous results on curvature-driven flows (Bloore 1977; Domokos 2014; Domokos and Gibbons 2012; Gage 1983; Grayson 1987) show that several of the mentioned three shape descriptors
evolve in a  monotonic or quasiconvex manner (for a partial overview on these results see Table 1 in the work of Miller et al. (2014)).
Our main goal in this paper is to establish analogous results
for the evolution of shape descriptors under distance-driven flows. Unless indicated
otherwise, we restrict ourselves to the description of the evolution of $C^3$-smooth, convex shapes either
in two dimensions (planar disks) or in three dimensions (solids).

In Subsection \ref{ss:models} we introduce curvature-driven and distance-driven flows
as models for collisional and frictional abrasion, respectively and also relate them to  the Aristotelian theory of abrasion.
In Subsection \ref{ss:descriptors} we list the aforementioned three types
of  shape descriptors and briefly review earlier results on their evolution under curvature-driven flows.

\subsection{Model Types: Curvature- and Distance-Driven Flows}\label{ss:models}

The mathematical framework for abrasion models are geometric partial differential equations (PDEs) describing the
evolution of shapes in time. One convenient way to write such a PDE is the so-called local notation where at each point
of the abrading surface the attrition speed $v$ in the direction of the inward surface normal is given. Alternatively,
in the global notation the abrading surface is identified by a scalar distance $r$ and the time derivative $r_t$
is expressed. (Here and henceforth subscripts refer to partial differentiation.) Based on these concepts
we identify two special classes of PDEs which appear to be particularly relevant as models of abrasion processes.
We call a PDE a curvature-driven flow if in the local notation it can be written in two and three dimensions, respectively, as
\begin{equation} \label{curv}
v=v(\kappa), \hspace{1cm} v=v(\kappa,\lambda),
\end{equation}
where in the two-dimensional case $\kappa$ is the scalar curvature, in the three-dimensional case $\kappa, \lambda$ are the principal curvatures.
Alternatively, we call a PDE a distance-driven flow if in the global notation it can be written as
\begin{equation} \label{dist}
r_t=f(r).
\end{equation}
Strictly speaking, Eq. (\ref{dist}) is not a partial differential equation, rather, a continuum of decoupled
ordinary differential equations. This evolution equation admits different interpretations:
if $r$ is interpreted as a radial distance measured from  a fixed point then we call Eq. (\ref{dist}) a radial distance-driven flow,
if $r$ is measured from a plane then we call it a parallel distance-driven flow (and in the latter case we often write Eq. (\ref{dist}) as $z_t=f(z)$).

Curvature- and distance-driven flows have been broadly investigated in the mathematical literature, here we only refer
to some fundamental papers. Distance-driven flows as models of particle abrasion have classical origins:
Aristotle postulated that particle abrasion is governed by a radial distance-driven flow of the type in Eq. (\ref{dist}),
in particular, he claimed (Krynine 1960) that if the function $f(r)$ is monotonically decreasing (with $f(r)<0$) then all shapes converge to the circle.
This model has never been verified from the physical point of view, however, we will show (Theorem~\ref{thm:limes}) that, if $df/dr=0$ at $r=0$ and $f(r)<0$, the mathematical claim is true, even for non-monotonic $f(r)$.
In a different physical context distance-driven flows also emerge in the so-called sharp interface limit (Kohn et al. 2006) of the Allen-Cahn equation, describing order-disorder transitions.

The study of curvature-driven flows was initiated much later by Lord Rayleigh (1942, 1944, 1944) who found that
ellipsoids are evolving in a self-similar fashion under the special curvature-driven flow
\begin{equation} \label{Rayleigh}
v=K^{\frac{1}{4}},
\end{equation}
where $K=\kappa\lambda$ is the Gaussian curvature. As in case of Aristotle's model, there is no physical argument behind Eq. (\ref{Rayleigh}),
however, the mathematical claim is undoubtedly correct. Firey (1974) proposed the first physically motivated curvature-driven model
for the abrasion of particles colliding in uniformly random directions with an infinite plane. Firey's model can be written as
\begin{equation} \label{Firey}
v=cK,
\end{equation}
(where $c$ is a constant) and under symmetry assumptions Firey proved that all convex shapes converge to the sphere under Eq. (\ref{Firey}). Andrews (1999) generalized Firey's argument to non-symmetrical shapes. The mathematical framework modelling abrasion by particles of arbitrary size was set up by Bloore (1977) who derived
the following curvature-driven flow
\begin{equation} \label{Bloore}
v=1+2bH+cK,
\end{equation}
where $b,c$ are constants and $H=0.5(\kappa + \lambda)$ is the mean curvature. The constants $b,c$ have been later identified in a paper of V\'arkonyi and Domokos (2011) based on results of Schneider and Weil (2008) as
\begin{equation}
b=\frac{M}{4\pi} \mbox{ and } c=\frac{A}{4\pi},
\end{equation}
where $M$ is the integrated mean curvature and $A$ is the surface area of the abrading particles. Note that Firey's model corresponds to the case of infinitely large abraders where the third term dominates in Bloore's model.

Bloore's model in Eq. (\ref{Bloore}) has also been studied by using a heuristic approximation by a system of ordinary differential equations, called the box equations (Domokos and Gibbons 2012, 2013), where Eq. (\ref{Bloore}) is reduced to the evolution of the principal axes of ellipsoids.
While the latter are certainly not invariant under the Eq. (\ref{Bloore}), this approximation still yields some qualitative insights, that is, in the box model it could be proven that all shapes converge to the sphere. Since Bloore's equation incorporates all collisional effects,
this result suggests that the frequently observed, non-spherical, elongated shapes of coastal pebbles are formed by a frictional process. The latter is non-local, as abrasion modes depend not just on local properties of the surface but also on global shape characteristics (Domokos and Gibbons 2012). From the mathematical point of view, the simplest non-local models are distance-driven flows of the type in Eq. (\ref{dist}).
In the paper of Domokos and Gibbons (2012) a set of axioms was proposed for such friction models and some were investigated in the frame of the box equations. In the current paper we aim to extend the analysis of some simple distance-driven models  to the full flow.

\subsection{Geological Shape Descriptors and Summary of Main Results}\label{ss:descriptors}

One central question in mathematical abrasion models is whether one can identify quantities which vary monotonically (or in some other, predictable manner) with time and
which can be measured in field campaigns or laboratory experiments. There are three types of such quantities which may be candidates for
both mathematical and experimental studies. We describe them below
and  we also summarize what is known about their evolution under curvature-driven flows (serving as models for collisional abrasion), then we add our new results which are related to some simple distance-driven flows
(serving as models corresponding to  the Aristotelian theory of abrasion and  as potential models for frictional abrasion).

\begin{itemize}

\item Axis ratios. In field and laboratory measurements, each particle is associated with three orthogonal axes  $a>b>c$ and the axis ratios
$y_1=c/a, y_2=b/a$ are regarded as geological shape descriptors. While the manual protocols for measuring these axes slightly vary,
according to all protocols the axes associated with a tri-axial ellipsoid are the actual principal axes of the ellipsoid.
In general, little is known about the evolution of axis ratios under curvature-driven flows of the type in Eq. (\ref{Bloore}), some partial understanding has been gained from the box equations (Domokos and Gibbons 2012, 2013) which can predict the  evolution of axis ratios for well-abraded, almost-ellipsoidal shapes. In this case it was found that
for collisional abrasion the evolution of axis ratios may be either monotonic with
$\lim_{t\to \infty}y_1=\lim_{t\to \infty}y_2=1$  or quasiconvex  with  $\lim_{t\to \infty}y_1=\lim_{t\to -\infty}y_1=\lim_{t\to \infty}y_2=\lim_{t\to -\infty}y_2=1$.
These results imply that, at least in this approximate model for well-abraded particles, all shapes ultimately converge to the sphere under collisional abrasion.
However, in the initial stages of abrasion where particles are still close to their fragmented original shape, very little is known on this subject.

 In our current paper we will only consider the  evolution of axis ratios if the initial shape is an ellipsoid. For radial distance flows of the type in Eq. (\ref{dist}), with $f(r)$ strictly convex or strictly concave,
in Sect. \ref{sec:radial} we will prove Theorem~\ref{thm:radial_axis_monotone} stating that  if the evolution starts from ellipsoids then axis ratios are monotonic functions of time.  Nevertheless, we show that even for the simplest nonlinear radial flow $f= \alpha r+ \beta r^2$, starting with a suitable ellipsoid, any pair of limits in the range $[0,1]$ can be achieved for the axis ratios (Proposition~\ref{prop:radial_axis_limits}).

In Sect.~\ref{sec:main} we consider orthogonal affinity as the simplest parallel distance-driven flow.  In Subsection~\ref{ss:dist1} we prove that the axis ratio of an ellipse evolves as a quasiconcave function  (Theorem~\ref{thm:parallel_planar_axes})
and point out that this function has a non-smooth maximum if the direction of the affinity is aligned with any of  the principal directions of the ellipse. In three dimensions we only consider the case
when one of the principal directions is orthogonal to the direction of the affinity. In this case we point out (Remark~\ref{ellipsoid_axis_ratio_affinity}) that the smaller axis ratio $y_1$ is a quasiconcave function, however, the larger axis ratio $y_2$
may have several extrema.

\item Isoperimetric ratio. In two dimensions, the isoperimetric ratio is defined as $I=(4\pi A)/P^2$ where $A$ is the enclosed area and $P$
is the perimeter of the curve. In three dimensions, we have  $I=(6\sqrt{\pi} V)/(A^{3/2})$ where $V$ is the volume and $A$ is the surface area of the solid.
The isoperimetric ratio $I$ has been measured both in the field (Miller et al. 2014) and in laboratory experiments (McCubbin et al. 2014).
The isoperimetric ratio  $I$ is of particular interest, because in case of the $v=\kappa$ curvature-driven flow (serving as a special model of collisional abrasion) it was proven by Gage (1983) that $I(t)$ is growing monotonically in time.

In case of radial distance flows, the results for axis ratios do not apply to the  isoperimetric ratio. As we point out in Remark~\ref{rem:radial_iso}, even in the case of planar, $D_2$-symmetric shapes and convex function $f(r)$, the evolution of $I(t)$ may be non-monotonic.

In case of parallel flows, in Sect.~\ref{sec:main}, Subsection~\ref{ss:dist2} we prove
for smooth convex bodies in arbitrary dimensions under orthogonal affinity that the isoperimetric ratio $I(t)$ is a quasiconcave function.
(Theorem~\ref{thm:quasiconvex}). In addition, we also show that there exists an orthogonal basis along the directions of which $dI/dt=0$ (Theorem~\ref{thm:basis}).

\item Number $T$ of static equilibrium points. Another, recently investigated indicator of abrasion processes is the number $T$ of static equilibrium points (Domokos et al. 2010). We assume that the abrading particle is a smooth convex body described by the scalar distance $r$ measured from the center of mass. In 2 dimensions, $r$ can be conveniently parametrized by the polar angle $\phi$ and in three dimensions by the Euler angles $(\phi,\theta)$, so the evolution of these shapes is given in two and three dimensions by the functions $r(\phi,t)$ and $r(\phi, \theta, t)$, respectively. Static equilibrium points are associated with spatial critical points of the aforementioned scalar functions, i.e. in two dimensions they are characterized by $r_{\phi}(\phi, t)=0$ and in three dimensions by $r_{\phi}(\phi, \theta,t)=r_{\theta}(\phi, \theta,t)=0.$  If $r$ is a Morse function (i.e. the shape of the particle is generic) then in two dimensions, (planar disks)
based on the sign of $r_{\phi\phi}$ we distinguish between stable
and unstable equilibria, and denote their numbers respectively by $S$ and $U$ and we have $S+U=T$. In three dimensions (solids), based on the eigenvalues of the Hessian of $r$,
 we distinguish between stable, unstable and saddle-type equilibria,
and denote their numbers respectively by $S,U,H$ and on generic (Morse) surfaces we have $S+U+H=T$. These numbers are related by the
the Poincar\'e-Hopf Theorem (Arnold, 1998) on topological invariants as
\begin{equation} \label{Poincare}
S-U=0  \mbox{  (in 2 dimensions)} ,\quad S+U-H=2 \mbox{  (in 3 dimensions)}.
\end{equation}
The pair $\{S,U\}$ is called the primary equilibrium class of the body while the Morse-Smale complex associated with the gradient of the distance function $r$ defines
the secondary and tertiary equilibrium classes of the body (Domokos et al. 2016a, 2016b).
The evolution of these numbers has been already measured in the field (Miller et al. 2014).
In case of the planar $v=\kappa$ flow (also called the curve shortening flow), it follows from  Grayson's result (Grayson 1987) that, if the reference point of $r(\phi,t)$ is fixed and it coincides with the center of mass, then $T(t)$
is monotonically decreasing. With some weakening assumptions on genericity and stochasticity, this
statement was generalized in a paper of Domokos (2014), showing that if $v_{\kappa}>0$ (two-dimension) and $v_{\kappa}, v_{\lambda}>0$ (three-dimension) then $T(t)$ can be approximated by a stochastic process the expected value of which is monotonically decreasing in time.

In  Sect.~\ref{sec:radial}  we will prove Theorem~\ref{thm:constantT} stating that the number $T$ of spatial critical points of $r$ evolving under Eq. (\ref{dist}) remains constant and this implies that if a convex body is evolving under the radial flow in Eq. (\ref{dist}) and the flow leaves the center of mass invariant then $T(t)=\mathrm{constant}$.
(The invariance of the center of mass is guaranteed by a sufficient symmetry group (e.g. $Z_2\times Z_2 \times Z_2$). If the center of mass is not invariant then its motion may be modeled as a white noise with zero expected value, added to $T(t)$ (Domokos 2014); in this case $T$ will be a random variable with constant expected value.)
Our argument also shows that all equilibrium classes (primary, secondary and tertiary) remain invariant under radial distance-driven flows if the center of mass does not move.

In Sect.~\ref{sec:main}, Subsection~\ref{ss:dist3}, we show that under  orthogonal affinity, as time $t$ tends to  infinity, the number $U$ of unstable points is approaching its minimal value $U=2$ (Theorem~\ref{thm:large}), also implying that
for sufficiently small  (positive) values of $t$ the number $S$ of stable points is approaching its minimal value $S=2$. This result suggests that $T(t)$ evolves as a quasiconcave function. This is, however, not true, one can easily find counterexamples. Nevertheless, in a weaker sense the statement can be still salvaged:
using a stochastic approach (similarly to the result of Domokos (2014)), in Subsection~\ref{ss:dist4} we show for planar rectangles that the probability $p(t)$ that a random truncation with a straight line
results in an increase of $T(t)$ is a quasiconcave function of $t$ (Theorem~\ref{thm:stochastic}) and the maximum of $p(t)$  coincides with the maximum of $I(t)$.
\end{itemize}

The paper is organized as follows:
In Sects.~\ref{sec:radial} and \ref{sec:main}  we present the results related to axis ratios, the isoperimetric ratio and the number of equilibria in this order in separate subsections. While the length of these subsections may differ substantially, this principle helps to organize the results.
In Sect.~\ref{sec:sum} we summarize our results and also formulate a conjecture for orthogonal affinity about the global connection between $I(t)$ and $T(t)$ in an averaged sense. Table 1 summarizes the structure of the paper, lists the most important references for curvature-driven flows as well as the main results obtained in the current work.

\begin{table}
\begin{center}
\label{t:1}

\begin{tabular}{||c||c|c||c|c||} \hline
Flow  driven by:  & \multicolumn{2}{|c||}{\textbf{Curvature}} & \multicolumn{2}{|c||}{\textbf{Distance}}\\
\hline
\hline
\textbf{Shape descriptor $\downarrow$} & Gauss & Bloore &  Radial  &  Parallel \\ \hline  \hline
Axis ratios & Domokos and & Domokos and & Subsect. \ref{sec:radial_axis}, & Subsect.~\ref{ss:dist1}, \\
 &  Gibbons 2012, & Gibbons 2012, & Theorems~\ref{thm:limes},\ref{thm:radial_axis_monotone} & Theorem~\ref{thm:parallel_planar_axes}, \\
 &  Firey 1974 & Bloore 1977 &  &  Remark~\ref{ellipsoid_axis_ratio_affinity} \\ \hline

Isoperimetric   &  Gage 1983 &          &   Subsect.~\ref{sec:radial_iso},  & Subsect.~\ref{ss:dist2},\\
       ratio       &            &          &   Remark~\ref{rem:radial_iso}     & Theorems~\ref{thm:quasiconvex}, \ref{thm:basis}        \\ \hline
Number of  & Grayson 1987  & Domokos 2014 & Subsect. \ref{sec:radial_num}, & Subsect. \ref{ss:dist34}, \\
 equilibria           &               &              & Theorem \ref{thm:constantT}    & Theorems \ref{thm:large}, \ref{thm:stochastic} \\
\hline
\end{tabular}
\caption{Structure of the paper, principal references and main results}
\end{center}
\end{table}

\section{Radial Distance Driven Flows}\label{sec:radial}

In this section we develop some results for radial distance-driven flows evolving under Eq. (\ref{dist}).
Before stating our results, we prove a lemma that we are going to use in many proofs later.

\begin{lemma}\label{lem:monotonicity}
Consider Eq. (\ref{dist}) as an ordinary differential equation with the unknown function $r(t)$. Assume that the function $f(r)$ is $C^1$-class differentiable. Let $r_1(t)$ and $r_2(t)$  be two arbitrary solutions of Eq. (\ref{dist}) satisfying $r_1(t) < r_2(t)$ for some $t \in \Re$. Then $r_1(t) < r_2(t)$ for every value $t \in \Re$.
\end{lemma}

\begin{proof}
By the Picard-Lindel\"of Theorem, for any initial value $r(0)=r_0$ there is a unique solution $r(t)$ satisfying Eq. (\ref{dist}) .
Clearly, if $f(r_0)=0$, then $r(t)=r_0$ is a  solution of Eq. (\ref{dist}).

Assume that $f(r_0) > 0$, and let $r_1 < r_0 < r_2$ be the roots of $f$ closest to $y_0$, if they exist.
Then $r(t)$ is a strictly increasing function satisfying $\lim\limits_{t \to \infty} r(t) = r_2$, and $\lim\limits_{t \to  -\infty} r(t) = r_1$.
Furthermore, the solution belonging to any initial condition $\bar{r}(t)=\bar{r}$, where $r_1 < \bar{r} < r_2$ can be written as $\bar{r}(t) = r(t-t_0)$ for some constant $t_0 \in \Re$. A similar consideration can be applied if $r_1$ or  $r_2$ do not exist, or if $f(r_0) < 0$.
Thus, every solution of the differential equation is strictly monotonic or constant, and thus, if $r_1(t)$ and $r_2(t)$ are solutions satisfying $r_1(t) < r_2(t)$ for some value of $t$, then $r_1(t) < r_2(t)$ for every value of $t$.
\end{proof}

\begin{remark}
The claim of Lemma~\ref{lem:monotonicity} can be also understood from the point of view of dynamical systems. The phase space of Eq. (\ref{dist}) is one-dimensional, so, regardless of the number and position of critical points, the ordering of any point set $r_i$ is invariant under the flow for any finite time $t$.
\end{remark}

\subsection{Axis Ratios}\label{sec:radial_axis}

Let $K$ be a convex body in $\Re^3$, with $Z_2\times Z_2 \times Z_2$-symmetry, the planes of reflection symmetry coinciding with the three coordinate planes. Then we call the ellipsoid with axes contained in the coordinate axes and containing the points of $\bd K$ on the coordinate axes on its boundary, the ellipsoid approximating $K$ and we identify the axis ratios $y_1,y_2$  of $K$ (as defined in Subsection~\ref{ss:descriptors}) with the axis ratios of its approximating ellipsoid. We can define similarly the axis ratio for plane convex bodies with $Z_2\times Z_2$-symmetry.

Recall Aristotle's claim that a pebble surface abrading according to Eq. (\ref{dist}), where $f$ is negative for $r > 0$ and $f(0)=0$, approaches a sphere for large values of $t$. We show that under an additional condition, the claim is true. Nevertheless, we will see that his claim is not true in general.

\begin{theorem}\label{thm:limes}
Let $f$ be a $C^{k+1}$-class function satisfying $f(0)=f'(0)= \ldots = f^{(k)}(0)=0$ and $f^{(k+1)}(0) \neq 0$. Let $r(\uu,t)$ be a solution of
Eq. (\ref{dist}) for $\uu \in \Sph^2$, satisfying the initial condition $r(\uu,0)=r_0(\uu)>0$. Furthermore, let $M(t) = \max \{ r(\uu,t):\uu \in \Sph^2 \}$, and assume that $f(r) < 0$ on $[0,M(0)]$. Then $\frac{r(\uu,t)}{M(t)}$ uniformly converges to $1$ on $\Sph^2$ as $t \to \infty$.
\end{theorem}

\begin{proof}
Note that as $f(r) < 0$ on $[0,M(0)]$, for every $\uu \in \Sph^2$, $r(u,t)$ is a strictly decreasing function of $t$ that tends to $0$ as $t \to \infty$. Since $f(r)$ is $C^{k+1}$-class, for any $r \in (0,M(0))$, $f(r) = \frac{1}{(k+1)!} f^{(k+1)}(\zeta_r) r^{k+1}$ for some $\zeta_r \in (0,r)$.
Applying this equality for sufficiently small values of $r$, the inequality $f(r) < 0$, where $r \in [0,M(0)]$, implies that $f^{(k+1)}(0) < 0$.
By the continuity of this function, it also follows that $\beta_1 \leq \frac{f^{(k+1)}(r)}{(k+1)!} \leq \beta_2 < 0$ in an interval $[0,r_0]$ for some suitable values $r_0 > 0$, and $\beta_1$ and $\beta_2$.
Let $t_0 \geq 0$ be any value such that $M(t_0) \leq r_0$.
Then, $\beta_1 r^{k+1} \leq f(r) \leq \beta_2 r^{k+1}$ on $[0,M(t_0)]$.

Solving the differential equation $r'_t(\uu,t)= f(r(\uu,t))$, we obtain $t-t_0 = \int\limits_{r(\uu,t_0)}^{r(\uu,t)} \frac{1}{f(s)} \, ds$.
Thus, it follows that $\int\limits_{r(\uu,t_0)}^{r(\uu,t)} \frac{1}{\beta_1 s^{k+1}} \, ds \leq t-t_0 \leq \int\limits_{r(\uu,t_0)}^{r(\uu,t)} \frac{1}{\beta_2 s^{k+1}} \, ds$, which yields the inequalities
\begin{equation}\label{eq:approx}
\frac{r(\uu,t_0)}{\sqrt[k]{1-k \beta_1 r^k(\uu,t_0) (t-t_0)}} \leq r(\uu,t) \leq \frac{r(\uu,t_0)}{\sqrt[k]{1-k \beta_2 r^k(\uu,t_0) (t-t_0)}} .
\end{equation}
 This immediately yields that
\begin{equation}\label{eq:estimate}
\frac{r(\uu,t_0)}{M(t_0)} \cdot \sqrt[k]{\frac{1-k \beta_2 M^k(t_0) (t-t_0)}{1-k \beta_1 r^k(u,t_0) (t-t_0)}} \leq \frac{r(\uu,t)}{M(t)} \leq 1,
\end{equation}
for any $t > t_0$ and $\uu \in \Sph^2$, where we note that the same inequality holds if we replace $t_0$ by any larger value.

Let $\varepsilon > 0$ be arbitrary, and let $\beta_1(t) = \frac{1}{(k+1)!} \min \{ f^{(k+1)}(r) : r \in [0,M(t)] \}$ and $\beta_2(t) = \frac{1}{(k+1)!} \max \{ f^{(k+1)}(r) : r \in [0,M(t)] \}$.
Clearly, $\lim\limits_{t \to \infty} \beta_1(t)= \lim\limits_{t \to \infty} \beta_2(t) = \frac{f^{(k+1)}(0)}{(k+1)!} < 0$.
Thus, there is some $\bar{t} \geq t_0$ such that $\sqrt[k]{\frac{\beta_2(\bar{t})}{\beta_1(\bar{t})}} > 1-\frac{\varepsilon}{2}$.
On the other hand, as the limit of the left-hand side of Eq. (\ref{eq:estimate}) is $\sqrt[k]{\frac{\beta_2(t_0)}{\beta_1(t_0)}}$, there is some $t_1 \geq \bar{t}$ such that for any $t > t_1$ and $\uu \in \Sph^2$, we have $1-\varepsilon < \frac{r(\uu,t)}{M(t)} \leq 1$. This proves the assertion and it also implies that the geometric shape described by $r(u)$ in a polar coordinate system uniformly converges to the sphere.
\end{proof}

\begin{remark}\label{rem:approximation}

If $f'(0) = \alpha < 0$, $f''(0) = \ldots = f^{(k)}(0) = 0$ and \linebreak  $\beta_1 \leq \frac{f^{(k+1)(r)}}{(k+1)!} \leq \beta_2 < 0$, then we can still approximate $r(\uu,t)$ by
\[
r(\uu,0) \cdot \sqrt[k]{\frac{\alpha}{(\alpha+\beta_1 r^k(\uu,0)) e^{-k\alpha t}-\beta_1 r^k(\uu,0)}} \leq r(\uu,t) \leq \]\linebreak
\[\leq r(\uu,0) \cdot \sqrt[k]{\frac{\alpha}{(\alpha+\beta_2 r^k(\uu,0)) e^{-k\alpha t}-\beta_2 r^k(\uu,0)}} .
\]
\end{remark}

\begin{remark}\label{rem:measurement}
Note that if $f^{(k+1)}(r)$ is close to a constant (or equivalently, if $M(0)$ is small), then $\beta_1 \approx \beta_2$, and thus,
$\lim\limits_{t \to \infty} \frac{r(\uu,t)}{M(t)} \approx \sqrt[k]{\frac{\alpha r^k(\uu,0) + \beta_1 r^k(\uu,0) M_0^k}{\alpha M_0^k + \beta_1 r^k(\uu,0) M_0^k}}$. Hence, the limit shape can be estimated by measuring $r(\uu,t)$ and $M(t)$ for large values of $t$. If $\beta_1 = \beta_2$, then the approximation becomes an equality, which shows that, despite Aristotle's claim, the limit shape can be different from a sphere. We will elaborate on this in Theorem \ref{thm:radial_axis_monotone} where we show that axis ratios may be  monotonically both  decreasing and increasing. In Proposition~\ref{prop:radial_axis_limits} we will show that even for quadratic $f(r)$, any axis ratio may be achieved as a limit.\end{remark}

\begin{remark}
Note that Theorem~\ref{thm:limes}, and Remarks~\ref{rem:approximation} and \ref{rem:measurement} remain valid if $r$ denotes distance from a plane or a line, and also for planar figures.
\end{remark}

\begin{theorem}\label{thm:radial_axis_monotone}
Let $E(0)$ be an ellipsoid, with its axes on the coordinate axes. Let $E(t)$ be the family of convex bodies generated by the evolution starting at $E(0)$ under the radial distance-driven flow in Eq. (\ref{dist}), where $f(r)$ is strictly decreasing and strictly convex/concave, respectively, for $r > 0$, and $f(0)=0$. Then, depending on the convex/concave property of $f(r)$, both axis ratios of $E(t)$ are monotonically decreasing/increasing functions of time, respectively.
\end{theorem}

\begin{proof}
For any $t > 0$, $E(t)$ is symmetric to any coordinate plane. Thus, the semi-axes of the approximating ellipsoid of $E(t)$ coincide with the radii of $E(t)$ in the direction of the coordinate axes. Let these radii be $c(t) \leq b(t) \leq a(t)$; note that by Lemma~\ref{lem:monotonicity}, if $c(0) \leq b(0) \leq a(0)$, then the same inequalities hold for any value of $t$.

We prove that if $f(r)$ is strictly concave for $r > 0$, then the axis ratios are increasing; we may apply the same argument if $f(r)$ is strictly convex. First, we show that $r \mapsto \frac{f(r)}{r}$ is a strictly decreasing function of $r$ for
 $r > 0$. Indeed, note that $\left( \frac{f(r)}{r} \right)' = \frac{f'(r) r - f(r)}{r^2}= \frac{1}{r} \left( f'(r) - \frac{f(r)}{r} \right)$. On the other hand, as $f(r)$ is strictly concave, it follows that $f'(r) < \frac{f(r)}{r}$ for any $r > 0$, and thus, $\left( \frac{f(r)}{r} \right)' < 0$ implies that $\frac{f(r)}{r}$ is strictly decreasing.

Consider now, say, the axis ratio $\frac{b(t)}{a(t)}$, and observe that both $a(t)$ and $b(t)$ are strictly decreasing functions of $t$, and satisfy $0 < b(t) \leq a(t)$ for all values of $t$.
Then $\left( \frac{b(t)}{a(t)} \right)' = \frac{b'(t) a(t) - a'(t) b(t)}{a^2(t)} = \frac{b(t)}{a(t)} \left( \frac{f(b(t))}{b(t)} - \frac{f(a(t))}{a(t)} \right)$. The inequalities $0 < b(t) \leq a(t)$ and the monotonicity of $\frac{f(r)}{r}$ yield that $\frac{f(b(t))}{b(t)} > \frac{f(a(t))}{a(t)}$, and hence, $\left( \frac{b(t)}{a(t)} \right)' > 0$, which implies that $\frac{b(t)}{a(t)}$ is an increasing function of time.
\end{proof}

Note that Theorem~\ref{thm:limes} implies that the axis ratios of any convex body, evolving under the flow in Eq. (\ref{dist}), with $f(0)=f'(0)=0$ and $f''(0) < 0$, tend to $1$ as $t \to \infty$. Nevertheless, any pair of `reasonable' pair can be obtained as a limit, if we drop the condition that $f'(0)=0$, even using a  simple quadratic function as $f(r)$.

\begin{proposition}\label{prop:radial_axis_limits}
Let $E_0 = E(0)$ be an ellipsoid, with its axes on the coordinate axes. Let $E(t)$ be the family of convex bodies evolving under the radial distance-driven flow in Eq. (\ref{dist}), with $f(r)=\alpha r+ \beta r^2$. Let the axis ratios of $E(t)$ be $0 < y_1(t) \leq y_2(t) \leq 1$. Then for any $\alpha < 0$, $\beta < 0$ and $0 < y_1 \leq y_2 \leq 1$ there is an ellipsoid $E_0$ such that $\lim\limits y_i(t) = y_i$ for $i=1,2$.
\end{proposition}

\begin{proof}
Let the semi-axes of $E_0$ be $0 < c \leq b \leq a$.
Then $y_1 = \frac{\alpha c + \beta ac}{\alpha a + \beta ac}$ and $y_2 = \frac{\alpha b + \beta ab}{\alpha a + \beta ab}$ (cf. Remark \ref{rem:measurement}).
Thus, by choosing $c = \frac{\alpha a y_1}{\alpha + \beta a (1- y_1)}$, $b = \frac{\alpha a y_2}{\alpha + \beta a (1- y_2)}$ and arbitrary $a$, the desired limit ratios can be achieved.
Clearly, if $a > 0$, then $0 < \frac{\alpha a y_1}{\alpha + \beta a (1- y_1)} \leq \frac{\alpha a y_2}{\alpha + \beta a (1- y_2)}$.
On the other hand, $\frac{\alpha a y_2}{\alpha + \beta a (1- y_2)} = \frac{a y_2}{1+\frac{\beta}{\alpha} a (1-y_2)} < a y_2 < a$.
\end{proof}

\subsection{Isoperimetric Ratio}\label{sec:radial_iso}
The following example shows that, unlike axis ratios (cf. Theorem~\ref{thm:radial_axis_monotone}), the isoperimetric ratio does not necessarily change monotonically under a radial distance driven flow.

\begin{remark}\label{rem:radial_iso}
Consider a square of side length two, centered at the origin, and replace two opposite edges of it by semicircles of unit radius. Let the obtained stadium-like convex region be $K_0$. Truncate $K_0$ by a circle $C_{\alpha}$ of radius $2 \cos \alpha$, where $0 \leq \alpha \leq \frac{\pi}{4}$, and denote the truncated figure by $K(\alpha)$.
An elementary consideration shows that $C_{\alpha}$ cuts off two arcs of the unit semicircles, each with central angle $4 \alpha$.
Then it is a matter of computation to show that the isoperimetric ratio of $K(\alpha)$ is
\[
\frac{I(K(\alpha))}{4 \pi} = \frac{4+\pi-4\alpha \cos(2 \alpha)-2 \sin(2 \alpha)}{\left( 4 + 2\pi - 8 \alpha + 8 \alpha \cos \alpha \right)^2},
\]
which is a convex function of $\alpha$, with its minimum attained at some $0 < \alpha_0 < \frac{\pi}{4}$.
First, imagine that $f(r) = \left\{ \begin{array}{l} 0, \hbox{ if } 0 \leq r \leq 2 \cos \alpha_0 \\ -\infty \hbox{ if } r > 2 \cos \alpha_0 \end{array} \right.$. Then, applying the `flow' in Eq. (\ref{dist}) to $K_0$ we obtain $K(\alpha_0)$, which has a smaller isoperimetric ratio. Clearly, we may replace $f(r)$ with a negative, concave, analytic function satisfying $f'(0)=0$ while still satisfying this property. On the other hand, by Theorem~\ref{thm:limes}, for large values of $t$, the shape of the figure obtained from $K_0$ is `almost' a circle, which has a larger isoperimetric ratio. Thus, the isoperimetric ratio of $K_0$, under this flow, is not a monotonic function of $t$.
\end{remark}

\subsection{Number of Equilibria}\label{sec:radial_num}

\begin{theorem} \label{thm:constantT}
The total number $T$ of spatial critical points (added number of local minima, maxima and saddles) of the function $r(\mathbf{u},t)$, $\mathbf{u} \in \Sph^{n-1}$  does not change in time under the flow in Eq. (\ref{dist}).
\end{theorem}

\begin{proof}
Let $\mathbf{u}^{\star}$ be, say, a local maximum of $r(\mathbf{u},0)$. Then $\mathbf{u}^{\star}$ has a neighborhood $V$ in $\Sph^n$ such that for any $\mathbf{u} \in V$, $r(\mathbf{u},0) \leq r(\mathbf{u}^{\star},0)$. Since $r(\mathbf{u},t)$ is a solution of Eq. (\ref{dist}), Lemma~\ref{lem:monotonicity} implies that for any $\mathbf{u} \in V$, and any value of $t$, we have $r(\mathbf{u},t) \leq
r(\mathbf{u}^{\star},t)$. Thus, $\mathbf{u}^{\star}$ is a local maximum of the surface $r(\mathbf{u},t)$ for any fixed value of $t$. Replacing the role of $0$ with any other value of $t$ we see that if $\mathbf{u}^{\star}$ is a local maximum of $r(\mathbf{u},t)$ for an arbitrary value of $t$, then it is also a local maximum of $r(\mathbf{u},0)$. Thus, the number of local maxima does not change in time. It can be shown similarly that the number of local minima does not change in time. The fact that the number of saddle points does not depend on $t$ follows from the Poincar\'e-Hopf Theorem (cf. Eq. \ref{Poincare}).
\end{proof}

\begin{corollary} \label{cor:constantT}
If the center of mass coincides with the origin then Theorem \ref{thm:constantT} implies that the total number $T$ of static equilibrium points is invariant under the radial flow in Eq. (\ref{dist}). Moreover, in this case
the proof of Theorem~\ref{thm:constantT} establishes that both the number $S$ of stable equilibria (sinks) and the number $U$ of unstable equilibria (sources) is constant and
this implies that the number $H$ of saddles also remains constant.
\end{corollary}

The proof of Theorem \ref{thm:constantT} also yields the following, more general statement:

\begin{corollary}\label{cor:MorseSmale}
Let $\mathcal{M}$ be the Morse-Smale complex on $\Sph^2$, defined by the gradient flow of the Euclidean distance function $r(\uu,t)$. Then $\mathcal{M}$ does not change in time under the flow in Eq. (\ref{dist}).
\end{corollary}

We remark that both Theorem~\ref{thm:constantT} and Corollary~\ref{cor:MorseSmale} remain valid if $r$ denotes distance from a plane or a line.
In the nondegenerate case, the invariance of the number of spatial critical points can be extended to a more general class of flows

\begin{equation}\label{dist1}
r_t=f(r,r_{u_i},t), \quad \mathbf{u} \in \Sph^{n-1}.
\end{equation}

\begin{theorem}\label{lem:constantT_generic}
If $r(\mathbf{u},t)$ is $C^2$-class and $f$ is $C^1$-class, the number $T$ of spatial critical points of the function $r(\mathbf{u},t)$, $\mathbf{u} \in \Sph^{n-1}$  does not change at generic bifurcations in time under the flow in Eq. (\ref{dist1}).
\end{theorem}

\begin{proof}
Without loss of generality, we may assume that $n=2$.

Generic saddle-node bifurcations of critical points (Arnold 1998) of $r(u,t)$ in time are characterized by
\begin{equation} \label{p1}
r_u=r_{uu}=0,
\end{equation}
and we also have
\begin{equation} \label{p2}
r_{uuu}, r_{ut} \not =0.
\end{equation}
Based on Eq. (\ref{dist1}) we can write
\begin{equation} \label{p3}
r_{ut}=r_{tu}=f_rr_u+f_{r_u}r_{uu},
\end{equation}
Equations~(\ref{p1}) and (\ref{p3}) yield $r_{ut}=0$, however, this contradicts Eq. (\ref{p2}) so we see
that generic bifurcations may not occur.
\end{proof}


\begin{remark}
Since we are primarily interested in the flows in Eq. (\ref{dist1}) as mathematical models of physical processes, the absence
of generic bifurcations in the model suggests that if a physical process is governed by Eq. (\ref{dist1}) then $T(t)$ will be constant.
\end{remark}


\section{Parallel Distance-Driven Flows: Orthogonal Affinity}\label{sec:main}
In this section we only consider one particular parallel distance driven flow: orthogonal affinity which is defined by

\begin{equation}\label{aff}
z_t=-z
\end{equation}
and we investigate the evolution of shape descriptors under Eq. (\ref{aff}).
\subsection{Axis Ratios}\label{ss:dist1}
In case of axis ratios, we only consider ellipses and ellipsoids (the set of which is invariant under Eq. (\ref{aff})), and first we prove
\begin{theorem} \label{thm:parallel_planar_axes}
In case of ellipses  evolving under Eq. (\ref{aff}), if the $z$ direction does not coincide with any of the principal directions then the  axis ratio $y(t)$ is a smooth quasiconcave function. This function reaches its global maximum at a point where the angle of the axes of the ellipse with the axis of the affinity is $\frac{\pi}{4}$.
If the $z$ direction coincides with any of the principal axes then $y(t)$ has a single, non-smooth maximum and it is smooth otherwise.
\end{theorem}

\begin{proof}
Let $E=E(1)$ be an ellipse with semi-axes $\lambda_1 > \lambda_2 > 0$ such that the angle of its axes with the $x$ coordinate axis are $0 \leq \alpha \leq \frac{\pi}{2}$ and $\alpha+\frac{\pi}{2}$. Let $h_t : \Re^2 \to \Re^2$ be the orthogonal affinity, with the $x$-axis as its axis, and ratio $t > 0$, and set $E(t)=h_t(E)$. Since the proof is straightforward if $\alpha=0$ or $\alpha=\frac{\pi}{2}$, we may assume that $0 < \alpha < \frac{\pi}{2}$. Then the quadratic form corresponding to $E(t)$ is
\[
\left[ \begin{array}{cc}
x & z \end{array} \right]
\left[
\begin{array}{cc}
1 & 0 \\ 0 & \frac{1}{t}
\end{array}
\right]
\left[ \begin{array}{cc}
\cos(\alpha) & \sin(\alpha) \\
-\sin(\alpha) & \cos(\alpha)
\end{array} \right]
\left[ \begin{array}{cc}
\frac{1}{\lambda_1^2} & 0 \\ 0 & \frac{1}{\lambda_2^2}
\end{array} \right]
\left[ \begin{array}{cc}
\cos(\alpha) & -\sin(\alpha) \\
\sin(\alpha) & \cos(\alpha)
\end{array} \right]
\left[
\begin{array}{cc}
1 & 0 \\ 0 & \frac{1}{t}
\end{array}
\right]
\left[ \begin{array}{c}
x \\ z \end{array}
\right]
\]
and thus, its matrix is
\[
\frac{1}{\lambda_1^2 \lambda_2^2} \left[
\begin{array}{cc}
\lambda_1^2 \sin^2 (\alpha) + \lambda_2^2 \cos^2 (\alpha) & \frac{\sin(\alpha) \cos(\alpha)}{t} \left( \lambda_1^2 - \lambda_2^2 \right) \\
\frac{\sin(\alpha) \cos(\alpha)}{t} \left( \lambda_1^2 - \lambda_2^2 \right) & \frac{1}{t^2} \left( \lambda_1^2 \cos^2 (\alpha) + \lambda_2^2 \sin^2 (\alpha) \right)
\end{array}
\right] .
\]
The axis ratio of $E(t)$ is the root of the ratio of the two eigenvalues of this matrix. Hence, denoting the two eigenvalues of this matrix by $0 < \Lambda_1(t) \leq \Lambda_2(t)$, to prove that the axis ratio of $E(t)$ is quasiconcave it suffices to prove that $\frac{\Lambda_1(t)}{\Lambda_2(t)}$ is quasiconcave.

Elementary calculations yield that
\begin{equation}\label{axisratiodiff}
\left( \frac{\Lambda_1(t)}{\Lambda_2(t)} \right)' =
\frac{2 t \lambda_1^2 \lambda_2^2 \left( \lambda_1^2 \cos^2 (\alpha) + \lambda_2^2 \sin^2 (\alpha) - t^2 \left( \lambda_1^2 \sin^2 (\alpha) + \lambda_2^2 \cos^2 (\alpha) \right) \right)}{M},
\end{equation}
where $M$ is positive for every value of $t$.
This quantity is positive for $0 < t < t_0$ and negative for $t > t_0$, where $t_0=\sqrt{\frac{\lambda_1^2 \cos^2 (\alpha) + \lambda_2^2 \sin^2 (\alpha)}{\lambda_1^2 \sin^2 (\alpha) + \lambda_2^2 \cos^2 (\alpha)}}$. This yields the first part of the theorem.

To prove the second part, we substitute $t=1$ into Eq. (\ref{axisratiodiff}), and examine for which values of $\alpha$ will this quantity be equal to zero. If $t=1$, then the numerator of the right-hand side of Eq. (\ref{axisratiodiff}) is $2 \lambda_1^2 \lambda_2^2 \left( \lambda_1^2 - \lambda_2^2 \right) \cos(2 \alpha)$, which is zero if, and only if $\alpha = \frac{\pi}{4}$.
\end{proof}

\begin{remark}\label{rem:axes_monotonicity}
A more elaborate computation yields that if $0 < \alpha < \frac{\pi}{2}$, the semiaxes $\lambda_1(t)=\frac{1}{\sqrt{\Lambda_1(t)}}$ and $\lambda_2(t)=\frac{1}{\sqrt{\Lambda_2(t)}}$ of the ellipse $E(t)$ are strictly increasing functions of $t$, satisfying $\lim\limits_{t \to 0+0} \lambda_2(t) = 0$, $\lim\limits_{t \to \infty} \lambda_2(t) = \frac{\lambda_1 \lambda_2}{\sqrt{\sin^2 \alpha \lambda_1^2 + \cos^2 \alpha \lambda_2^2}}$, $\lim\limits_{t \to 0+0} \lambda_1(t) = \sqrt{\cos^2 \alpha \lambda_1^2 + \sin^2 \alpha \lambda_2^2}$ and $\lim\limits_{t \to \infty} \lambda_1(t) = \infty$. Furthermore, we note that if $\lambda_1 \neq \lambda_2$, then $\frac{\lambda_1 \lambda_2}{\sqrt{\sin^2 \alpha \lambda_1^2 + \cos^2 \alpha \lambda_2^2}} < \sqrt{\cos^2 \alpha \lambda_1^2 + \sin^2 \alpha \lambda_2^2}$. Observe that the first (smaller) constant is equal to  half of the length  of the interval of the $x$ axis inside the ellipse while the second, larger constant is equal to  half of the projection of the ellipse onto the $x$ axis.
\end{remark}

Even though the approach applied in the proof of Theorem~\ref{thm:parallel_planar_axes} works in any dimension, we could not modify it even for the $3$-dimensional case, due to computational difficulties. Nevertheless, following the ideas in the proof of Theorem~\ref{thm:parallel_planar_axes}, and using Remark~\ref{rem:axes_monotonicity}, we can prove the following.

\begin{remark} \label{ellipsoid_axis_ratio_affinity}
Let $E$ be an ellipsoid with semiaxes $A_1$, $A_2$ and $A_3$ of lengths $\lambda_1$, $\lambda_2$ and $\lambda_3$, respectively.
Assume that $A_1$ lies in the $x$ direction, and that the angle between $A_2$ and the $y$ direction is $0 \leq \alpha \leq \frac{\pi}{2}$.
Let $E(t)$ be the family of ellipsoids, evolving under Eq. (\ref{aff}), satisfying $E(1)=E$. Let $y_1(t) \leq y_2(t)$ be the axis ratios of $E(t)$, and set
\[
\bar{\lambda} =\sqrt{\cos^2 \alpha \lambda_2^2 + \sin^2 \alpha \lambda_3^2 - \sin \alpha \cos \alpha | \lambda_2^2-\lambda_3^2 |
\sqrt{\frac{\cos^2 \alpha \lambda_2^2 + \sin^2 \alpha \lambda_3^2}{\sin^2 \alpha \lambda_2^2 + \cos^2 \alpha \lambda_3^2}}} .
\]
Then $y_1(t)$ is a quasiconcave function of $t$. Furthermore, $y_2(t)$ is a quasiconcave function if, and only if
$\alpha = 0$, or $\alpha = \frac{\pi}{2}$, or $\lambda_2 = \lambda_3$, or $\lambda_1 \geq \bar{\lambda}$.
\end{remark}


\subsection{Isoperimetric Ratio}\label{ss:dist2}

Following Pisanski (1997), we call the quantity  $I(K) = \frac{\vol(K)}{\left( \surf(K) \right)^{\frac{n}{n-1}}} \cdot \frac{\left( \surf(\B^n) \right)^{\frac{n}{n-1}}}{\vol(\B^n)}$  the isoperimetric ratio of the $n$-dimensional convex body $K$, where $\B^n$ denotes the Euclidean unit ball with the origin as its center. (Note that for $n=2,3$ this definition yields the formula provided in Subsection~\ref{ss:descriptors}.) We  remark that other variants of this concept are also used in the literature (e.g., Ball 1991; Firey 1960; Green 1953).
We denote by $H$ a hyperplane passing through the center of mass of $K$, with normal vector $\vv \in \Sph^{n-1}$,
let $a^v_t : \Re^n \to \Re^n$ be the orthogonal affinity, with $H$ as its fixed hyperplane, and ratio $t > 0$, and set $K^v(t)= a^v_t(K)$.

\begin{theorem}\label{thm:quasiconvex}
For every convex body $K \subseteq \Re^n$ and every $\vv \in \Sph^{n-1}$, the isoperimetric ratio of $K^v(t)$ is a quasiconcave function of $t$.
\end{theorem}

\begin{proof}
Set $V(t) = \vol K^v(t)$ and $A(t) = \surf K^v(t)$.
Without loss of generality, we may assume that $V(1) = 1$ which implies that $V(t) = t$.
Thus we need to show that the function $I(t) = \frac{t}{A^s(t)}$, where $s=\frac{n}{n-1}$, is quasiconcave.
We show that $I'(t)$ has exactly one root, which, since $I(t) \to 0$ as $t \to 0$ or $t \to \infty$, yields that here $I(t)$ has a maximum.

Observe that $I'(t) = \frac{A(t) - stA'(t)}{A^{s+1}(t)}$, and that $\left( A(t) - stA'(t) \right)'= (1-s)A'(t) - stA''(t)$.
Note that as $1-s<0$ and $A'>0$, if $A''(t)>0$, then the numerator of $I'(t)$ is strictly decreasing, which, combined with the limits
$I(t) \to 0$ as $t \to 0$ or $t \to \infty$, implies the assertion.

We show that $A''(t)>0$.
Let us imagine $H$ as the hyperplane $\{ x_n = 0 \}$, let $h: \Re^n \to H$ be the orthogonal projection onto $H$, and set $h(K)=D$.
Then $\bd K$ is the union of the graphs of two functions defined on $D$, and the set $K \cap h^{-1}(\bd D)$.
Let these two functions be $x_n = f(x_1,x_2,\ldots,x_{n-1})$ and $x_n = g(x_1,x_2,\ldots,x_{n-1})$.
Then
\[
A(t) = \int_{x \in D} \sqrt{1+t^2 \grad^2 f} \dif x + \int_{x \in D} \sqrt{1+t^2 \grad^2 g} \dif x + \int_{x \in \bd D} t \left( h^{-1}(x) \cap K \right) \dif x.
\]

Observe that $\frac{\dif^2}{\dif t^2} \int_{x \in \bd D} t \left( h^{-1}(x) \cap K \right) \dif x= \int_{x \in \bd D} \left( h^{-1}(x) \cap K \right) \dif x \frac{\dif^2 t}{\dif t^2} = 0$.
Furthermore,
\[
\frac{\dif^2}{\dif t^2} \int_{x \in D} \sqrt{1+t^2 \grad^2 f} \dif x = \int_{x \in D} \frac{\dif^2}{\dif t^2} \sqrt{1+t^2 \grad^2 f} \dif x =
\]
\[
= \int_{x \in D} \frac{\grad^2 f}{\left( 1+t^2  \grad^2 f \right)^{3/2}} \dif x,
\]
where the integrand is positive, which implies that the integral is also positive. We obtain similarly that $\frac{\dif^2}{\dif t^2} \int_{x \in D} \sqrt{1+t^2 \grad^2 g} \dif x > 0$. This yields the assertion.
\end{proof}

\begin{theorem} \label{thm:basis}
Let $K \subset \Re^n$ be a convex body. For any $\vv \in \Sph^{n-1}$, let $I_v(t)$ denote the isoperimetric ratio of $a^v(K(t))$. Then there is an orthonormal basis $\ee_1, \ee_2, \ldots, \ee_n \in \Re^n$ such that for any $i$, $I'_{e_i}(1) = 0$.
\end{theorem}

\begin{proof}
For any $\vv \in \Sph^{n-1}$ let $f(\vv) = I'_v(1)$. Then, since $\bd K$ is $C^2$-class differentiable, $f : \Sph^{n-1} : \Re$ is continuous.
We show that for any orthonormal basis $\ee_1, \ee_2, \ldots, \ee_n$, we have $\sum_{i=1}^n f(\ee_i) = 0$.

Given an orthonormal basis, we define a function $g: \Re^n \to \Re$ in the following way: For $i=1,2,\ldots, n$, let $H_i$ be a hyperplane orthogonal to $\ee_i$, and
$a^i_t$ the orthogonal affinity, of ratio $t$, with $H_i$ as its fixed hyperplane.
We set $K(t_1,t_2,\ldots,t_n) = a^1_{t_1}(a^2_{t_2}(\ldots a^n_{t_n}(K) \ldots))$, and note that this body is independent of the order in which the affinities are carried out. Finally, let $g(t_1,t_2\ldots,t_n)$ be the isoperimetric ratio of the body $K(t_1,t_2,\ldots,t_n)$. Clearly, this function is differentiable at $(1,1,\ldots,1)$. Let $\vv = \sum_{i=1}^n \ee_i$. Then, by the linearity of the directional derivatives, we have $g'_v(1,\ldots,1) = \sum_{i=1}^n g'_{e_i}(1,\ldots,1) = \sum_{i=1}^n f(e_i)$. On the other hand, $g(t,\ldots,t)=g(1,\ldots,1)$ clearly holds for every $t > 0$, which yields that $g'_v(1,\ldots,1)=0$.

Now we prove the following, more general statement: If $h : \Sph^{n-1} : \Re$ is a continuous function such that for every orthonormal basis
$\ee_1, \ee_2, \ldots, \ee_n$, we have $\sum_{i=1}^n h(\ee_i) = 0$, then there is an orthonormal basis $\ee_1, \ee_2, \ldots, \ee_n$ such that $h(\ee_i)=0$ for every value of $i$.
This clearly implies the assertion.

We show this statement by induction on $n$. If $n=1$, then the statement is obvious. Now, assume that the statement holds for functions defined on $\Sph^{n-2}$.
Consider some $h : \Sph^{n-1} : \Re$ satisfying our conditions. Then there are some (orthogonal) vectors $\uu,\mathbf{w} \in \Sph^{n-1}$ such that $h(\uu) \leq 0 \leq h(\mathbf{w})$.
Thus, there is some  $\vv \in \Sph^{n-1}$ such that $h(\vv)=0$. We identify the set of vectors in $\Re^n$, perpendicular to $\vv$ with the space $\Re^{n-1}$, and let $h_v$ denote restriction of $h$ to this subspace.
Then, if $\ee_1,\ldots,\ee_{n-1} \in \Sph^{n-2}=\Re^{n-1} \cap \Sph^{n-1}$ is an orthogonal basis in $\Re^{n-1}$, then, adding $\vv$ to it we obtain an orthogonal basis in $\Re^n$, and thus, we obtain that $\sum_{i=1}^{n-1} h(\ee_i) = h(\vv) + \sum_{i=1}^{n-1} h(\ee_i) = 0$. Hence, we can apply the inductive hypothesis to $h_v$, which yields the required statement.
\end{proof}

Note that we have proven the following, stronger statement.

\begin{corollary}
Any orthonormal $k$-frame $\ee_1, \ee_2, \ldots, \ee_k$ in $\Re^n$, satisfying $I'_{\ee_i}(1) = 0$ for every value of $i$, can be completed to an orthonormal basis $\ee_1, \ldots,\ee_n$ satisfying the same property.
\end{corollary}

Theorem~\ref{thm:basis} suggests that if a plane convex body $K$ is symmetric to two perpendicular axes then, using the notations of Theorem~\ref{thm:basis}, $I'_v(1) = 0$ if the angle of $v$ and the symmetry axes of $K$ is $\frac{\pi}{4}$.

\begin{remark}\label{rem:45degrees}
Let $K$ be a plane convex body symmetric to the line $y=x$, and let $K(t)$ denote the image of $K$ under the orthogonal affinity defined by $(x,y) \mapsto (x,ty)$. Then $\left. \frac{\dif}{\dif t} I(K(t)) \right|_{t=1} = 0$.
\end{remark}

\begin{proof}
For simplicity, let $I(t)=I(K(t))$, and assume that $\area(K) = 1$. Then, by the proof of Theorem~\ref{thm:quasiconvex}, $I'(t) = \frac{P(t) - 2tP'(t)}{P^{2}(t)}$, where $P(t)$ is the perimeter of $K(t)$.

Let $\bd K$ be defined by the polar curve $\phi \mapsto r(\phi)$, where $0 \leq \phi \leq 2\pi$. Then the symmetry of $K$ implies that $r(\phi)=r\left( \frac{\pi}{2} - \phi \right)$ for every value of $\phi$.
The parametric form of $\bd(K(t))$ is $(r(\phi) \cos (\phi), t r(\phi) \sin \phi)$, $0 \leq \phi \leq 2\pi$.
Thus,
\[
P(t) = \int\limits_0^{2\pi} \sqrt{(r'(\phi) \cos \phi - r(\phi) \sin \phi)^2 + t^2 (r'(\phi) \sin \phi + r(\phi) \cos \phi)^2} \, d \phi,
\]
and
\begin{equation}\label{45_1}
P'(1) =
\int\limits_0^{2\pi} \frac{(r'(\phi) \sin \phi + r(\phi) \cos \phi)^2}{\sqrt{r(\phi)^2 + (r'(\phi))^2}} \, d \phi
\end{equation}
Using the substitution $\phi = \frac{\pi}{2} - u$ and the identities $r(\phi)=r\left( \frac{\pi}{2} - \phi \right)$ and
$r'(\phi)=-r'\left( \frac{\pi}{2} - \phi \right)$, this yields that
\begin{equation}\label{45_2}
P'(1) =
\int\limits_0^{2\pi}
\frac{(-r'(\phi) \cos \phi + r(\phi) \sin \phi)^2}{\sqrt{r(\phi)^2 + (r'(\phi))^2}} \, d \phi
\end{equation}

Combining Eqs. (\ref{45_1}) and (\ref{45_2}), we obtain that
\[
2 P'(1) =
\int\limits_0^{2\pi} \sqrt{r(\phi)^2 + (r'(\phi))^2} \, d \phi = P(1).
\]
This, together with $I'(1) = \frac{P(1) - 2P'(1)}{P^{2}(1)}$, yields the assertion.
\end{proof}

\subsection{Number of Equilibria} \label{ss:dist34}
\subsubsection{Deterministic Results}\label{ss:dist3}

Throughout this subsection, we assume that $K \subset \Re^n$ is a convex body with nowhere vanishing curvature. Note that this condition implies, in particular, that $K$ is strictly convex.

\begin{theorem}\label{thm:large}
Let $U_K(t)$ denote the number of unstable points of $K(t)$, with respect to its center of mass.
If $t$ is sufficiently large, then $U_K(t) = 2$.
\end{theorem}

\begin{proof}
First, we prove the assertion for $n=2$.

Let the origin $o$ be the center of mass of $K$, and $H$  be the $x$-axis.
Let the projection of $K$ onto the $x$ axis be $[a,b]$, with $a < 0$ and $b > 0$.
Then $\bd K$ is the union of the graphs of two functions, defined on $[a,b]$, one strictly concave and the other one strictly convex.
Let $f : [a,b] \to \Re$ be the strictly concave one, which is then $C^2$-class, and has nonvanishing curvature.
Then $\bd K(t)$ can be written as the set $\{ (x,tf(x)) \in \Re^2: x \in [a,b]$.

Let $x_0 \in (a,b)$ be the value with $f'(x_0) = 0$, and note that $x_0$ uniquely exists, as under our conditions, $f'$ is strictly decreasing.
First, observe that $K(t)$ has an equilibrium at $\left( x,tf(x) \right)$ if, and only if the tangent line of $f$
at this point is perpendicular to the position vector of the point, or in other words, if
\begin{equation}\label{eq:planar}
0 = \langle (x,tf(x)), (1,tf'(x))\rangle = x+t^2 f(x) f'(x) .
\end{equation}
For any fixed $x$ with $f'(x) \neq 0 \neq f(x)$, the right-hand side expression is a strictly monotonous function of $t$, which means that it is satisfied for exactly one value of $t$.
Thus for any $\varepsilon > 0$, if $t$ is sufficiently large, then we have one of the following for any equilibrium point $(x,tf(x))$ of $K(t)$ on $f$:
\begin{enumerate}
\item $|x - x_0 | < \varepsilon$;
\item $|x-x'_0| < \varepsilon$ for some $x'_0$ satisfying $f(x'_0)=0$.
\end{enumerate}

First, we consider the first type equilibria.
Note that $f(x_0) >0$, $f'(x_0)$ and by concavity and the nonvanishing of curvature, $f''(x_0) < 0$.
Thus, there is some value $t_0$ such that $t_0^2 \left( (f'(x_0))^2 + f(x_0) f''(x_0) \right) < -1$.
Furthermore, since $f(x) f''(x) + (f'(x))^2$ is a continuous function of $x$,
there is some neighborhood $V$ of $x_0$ such that $t_0^2 \left( (f'(x))^2 + f(x) f''(x) \right) < -1$ holds for any $x \in V$.
Note that for any $x \in V$, the same inequality holds for any $t > t_0$.

As the $x$-derivative of the right-hand side $R_t(x)$ of Eq. (\ref{eq:planar}) is \linebreak
$1 + t^2 \left( (f'(x))^2 + f(x) f''(x) \right) $, we obtain that for any $t \geq t_0$, $R_t(x)$
is a strictly decreasing function on $V$.
Hence, choosing $\varepsilon> 0$ such that $(x_0-\varepsilon,x_0+\varepsilon) \subset V$, for every sufficiently large $t \geq t_0$,
Eq. (\ref{eq:planar}) is satisfied for exactly one value $x_t$ of $x \in (a,b)$. Furthermore, in this case for any $x \in U$, $x < x_t$
yields $R_t(x) > 0$ and $x > x_t$ yields $R_t(x) < 0$,  from which it follows that $tK$ has an unstable point at $(x_t,tf(x_t))$.

Now we consider a second type equilibrium point $x'_t$.
Applying a similar argument as in the previous case, one can see that if $t$ is sufficiently large, then, in a neighborhood of $x'_0$,
the Euclidean distance function of $\bd K$ is minimal at $x_t$, which yields that $t$ is a stable point of $K$.
Thus, if $t$ is sufficiently large, then $U_K(t) = 2$.

Now, we prove the statement for any dimension $n > 2$.
Let $H = \{ x_n = 0 \}$, which we identify with $\Re^{n-1}$, and let $K_0$ be the orthogonal projection of $K$ onto $H$.
Let $f : K_0 \to \Re$ be the strictly concave function defining ``one half'' of $\bd K$.
Note that at any point of $f$, the supporting hyperplane of $K$ is spanned by the vectors $(0,\ldots, 0,1,0, \ldots, 0, \ldots, t\partial_i f)$,
where the $i$th coordinate is $1$, and $i=1,2,\ldots, n-1$.
Thus, the equilibria correspond to the points where the vector $(x_1,x_2,\ldots,x_{n-1},tf)$ is perpendicular to each of these vectors, and hence
to the solutions of the system of equations
\[
x_i + t^2 f \partial_i f = 0.
\]
From now on we need to repeat the planar argument.
\end{proof}

A more elaborate version of our argument yields the following, stronger statement in higher dimensions. To formulate it, let $eq_i(L)$ denote the number of equilibria of the convex body $L$, with exactly $i$ negative eigenvalues.
Note that  $eq_0(L)$ is the number of stable points, and $eq_n(L)$ is the number of unstable points.

\begin{corollary}\label{cor:allindices}
$K_0 = H \cap K$. Then, if $t$ is sufficiently large, $eq_n (K(t)) = 2$, and for any $i=0,1,\ldots, n-1$, $eq_i(K(t)) = eq_i(K_0)$.
\end{corollary}

\begin{remark}
We note that, following the idea of the proof of Theorem~\ref{thm:large}, one can obtain analogous statements to Theorem~\ref{thm:large} and Corollary~\ref{cor:allindices} about
the number $S_K(t)$ of stable equilibria in the case when $t> 0$ is sufficiently small.
\end{remark}

\subsubsection{Stochastic Results}\label{ss:dist4}


In this subsection we will first show that the probability that a randomly picked side of a polygon contained in a rectangle with sides  $2$ and $2a$ carries an equilibrium point is a quasiconcave function of $a$ and has its maximum exactly at $u=v$.

\begin{theorem}\label{thm:stochastic}
Let $R$ be a rectangle of side lengths $2$ and $2a$, where $a > 0$. Using uniform distribution,  choose two points $p_1,p_2 \in K$ independently. Let $p(a)$ denote the probability that the projection of the center of $R$ on the line containing $[p_1,p_2]$ is contained in the segment $[p_1,p_2]$. Then $p(a)$ is a quasiconcave function of $a$, which is maximal if, and only if $a=1$.
\end{theorem}

\begin{proof}
Since $p(a)=p \left( \frac{1}{a} \right)$, we may assume that $0 < a \leq 1$.
To make our computations simpler, we assume that the vertices of the rectangle are $(\pm 1, \pm a)$.

Let the two points be $p_1=(x_1,y_1)$ and $p_2=(x_2,y_2)$, with $x_1,x_2 \in [-1,1]$ and $y_1,y_2 \in [-a,a]$ be chosen uniformly.
Without loss of generality, we may assume that the slope of the line $L$ containing $p_1$ and $p_2$ is non-positive, and that
$L$ intersects the $y$-axis above $o$. Let $q$ denote the orthogonal projection of $o$ onto $L$.
Instead of the Cartesian coordinates, we use the following coordinate system: $r=|q|$, $m \geq 0$ is the slope of the line through $o$ and $q$,
$d_1$ is the signed distance of $p_1$ and $q$, and $d_2$ is the signed distance of $p_2$ and $q$, where the orientation of $L$ is chosen such that the signed distance of the intersection of $L$ with the line $y=a$ is positive from its intersection with the line $x=1$. Then $p_i=\left( \frac{r-d_i m}{\sqrt{1+m^2}}, \frac{m+d_i}{\sqrt{1+m^2}} \right)$ for $i=1,2$. The Jacobian of this coordinate transformation is $|J| = \frac{|d_2-d_1|}{1+m^2}$.

First, consider the case that $L$  separates $(1,a)$ from both $(-1,a)$ and $(1,-a)$.
Let $q$ be the projection of $o$ onto $L$.
An elementary computation yields that the conditions that $L$ separates $(1,a)$ from both $(-1,a)$ and $(1,-a)$, $d_2 \geq d_1$, and
that $q \in [p_1,p_2]$ are equivalent to the inequalities $m \geq 0$, $\frac{|1-am}{\sqrt{1+m^2}} \leq r \leq \frac{1+am}{\sqrt{1+m^2}}$, $\frac{r-\sqrt{1+m^2}}{m} \leq d_1 \leq a \sqrt{1+m^2}-rm$, $d_1 \leq d_2 \leq a \sqrt{1+m^2}-rm$.
Let $A= \int\limits_0^{a \sqrt{1+m^2}-rm} \int\limits_{\frac{r-\sqrt{1+m^2}}{m}}^0 \frac{d_2-d_1}{1+m^2} \dif d_1 \dif d_2$.
Then an elementary computation shows that if $0 < a \leq \frac{1}{\sqrt{8}}$, then the required probability is
\[
X_1 =  \frac{1}{2a^2} \left( \int\limits_0^a \int\limits_{\frac{|1-am|}{\sqrt{1+m^2}}}^{\sqrt{1+m^2}} A \dif r \dif m + \int\limits_a^{\frac{1-\sqrt{1-8a^2}}{4a}} \int\limits_{\frac{|1-am|}{\sqrt{1+m^2}}}^{\frac{a\sqrt{1+m^2}}{m}} A \dif r \dif m + \int\limits_{\frac{1+\sqrt{1-8a^2}}{4a}}^{\infty} \int\limits_{\frac{|1-am|}{\sqrt{1+m^2}}}^{\frac{a\sqrt{1+m^2}}{m}} A \dif r \dif m \right) ,
\]
and if $\frac{1}{\sqrt{8}}\leq a \leq 1$, then it is
\[
X_2 =  \frac{1}{2a^2} \left( \int\limits_0^a \int\limits_{\frac{|1-am|}{\sqrt{1+m^2}}}^{\sqrt{1+m^2}} A \dif r \dif m +
\int\limits_a^{\infty} \int\limits_{\frac{|1-am|}{\sqrt{1+m^2}}}^{\frac{a\sqrt{1+m^2}}{m}} A \dif r \dif m \right) .
\]

Next, we examine the case that $L$ separates the points $(1,\pm a)$ from the points $(-1,\pm a)$.
A similar computation shows that in this case the required probability is
\[
Y= \frac{1}{2a^2} \int_{\frac{1}{a}}^{\infty} \int_0^{\frac{am-1}{\sqrt{m^2+1}}}  \int_{\frac{r}{m}-\frac{m}{\sqrt{m^2+1}}}^0 \int_0^{\frac{r}{m}+\frac{m}{\sqrt{m^2+1}}} \frac{d_2-d_1}{1+m^2} \dif d_2 \dif d_1 \dif r \dif m.
\]

Now, assume that $L$ separates the points $(\pm 1,- a)$ from the points $(\pm 1, a)$, and set $B=\int\limits_{-\frac{a}{\sqrt{1+m^2}}-r}^0 \int\limits_0^{\frac{a}{\sqrt{1+m^2}}-r} \dif d_2 \dif d_1$.
Then the probability that $q \in [p_1,p_2]$ is
\[
Z=  \frac{1}{2a^2} \left( \int_0^{\frac{1}{a}-1} \int_0^{\frac{a}{\sqrt{1+m^2}}} B \dif r \dif m + \int_{\frac{1}{a}-1}^{\frac{1}{a}} \int_0^{\frac{1-am}{\sqrt{m^2+1}}} B \dif r \dif m \right).
\]

Finally, observe that since the slope of $L$ is nonnegative and $L$ intersects the $y$-axis above $o$, it does not separate the point $(-1,-a)$ from the other three vertices.
Thus
\[
p(a)= \left\{ \begin{array}{l} X_1+Y+Z, \hbox{ if } 0 < a \leq \frac{1}{\sqrt{8}},\\
X_2+Y+Z \hbox{ if } \frac{1}{\sqrt{8}}\leq a \leq 1 . \end{array} \right.
\]
Evaluating these integrals, numeric computations yield the assertion.
\end{proof}

\begin{figure}[ht]
\includegraphics[width=0.5\textwidth]{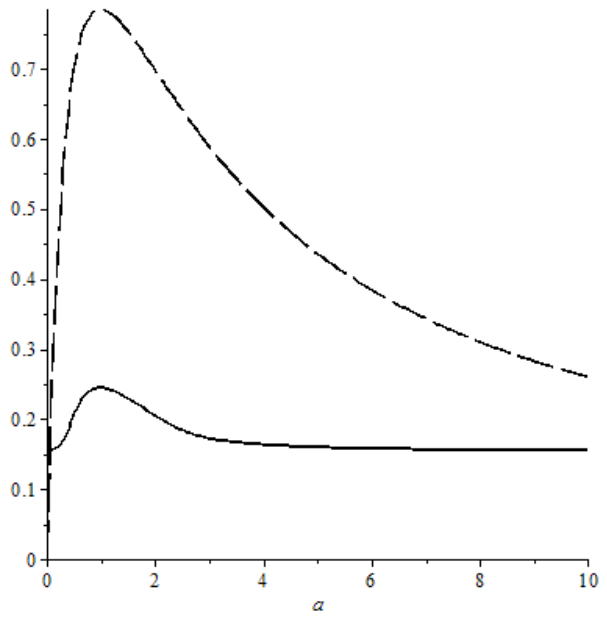}
\caption[]{Continuous line: the probability that a randomly chosen segment contains an equilibrium point, as a function of $a$. Dashed line: the isoperimetric ratio of the rectangle.}
\label{fig:stochastic}
\end{figure}

Figure~\ref{fig:stochastic} shows the value of $p(a)$, and also the isoperimetric ratio of the rectangle on the interval $(0,10]$.
Note that both functions attain their maxima at $a=1$.

\begin{remark}\label{rem:ellipse_stochastic}
We may examine the problem in Theorem~\ref{thm:stochastic} not only for rectangles but also for ellipses. More specifically, let $E$ be an ellipse of semi-axes $1$ and $a > 1$. Let us choose two points $p_1$ and $p_2$ randomly and independently in $E$, using uniform distribution.
Let $p(a)$ denote the probability that the orthogonal projection of the center of $E$ onto the line of $p_1$ and $p_2$ lies on the segment $[p_1,p_2]$. Then a similar computation to that in the proof of Theorem~\ref{thm:stochastic} yields that if $a > 1$,
\[
p(a) = \frac{4}{a^2 \pi^2} \int\limits_0^{\infty} \int\limits_0^{\frac{a\sqrt{1+m^2}}{\sqrt{1+a^2m^2}}} \int\limits_0^{D_2} \int\limits_{D_1}^0 \frac{d_2-d_1}{1+m^2} \, d d_1 \, d d_2 \, d r \, d m,
\]
where $D_1 = \frac{-(a^2-1)rm^2-am \sqrt{m^2+1} \sqrt{a^2+m^2-r^2(m^2+1)}}{a^2+m^2}$ and \linebreak $D_2 = \frac{-(a^2-1)rm^2+am \sqrt{m^2+1} \sqrt{a^2+m^2-r^2(m^2+1)}}{a^2+m^2}$. Nevertheless, due to computational difficulties, we could not prove a statement similar to Theorem~\ref{thm:stochastic}.
\end{remark}

\begin{question}\label{ques:ellipse_stochastic}
Can Theorem~\ref{thm:stochastic} be modified for ellipses instead of rectangles?
\end{question}

\begin{question}\label{ques:convexbody_stochastic}
Let $K$ be an origin-symmetric plane convex body $K$ and $L$ be a line through the origin. Let $h_a$ denote the orthogonal affinity with axis $L$ and ratio $a$. Define $p(a)$ similarly to that in Theorem~\ref{thm:stochastic}. Prove or disprove that for suitably chosen $K$ and $L$, $p(a)$ and the isoperimetric ratio of $h_a(K)$ attain their maxima at different values of $a$.
\end{question}

The problem in Theorem~\ref{thm:stochastic} can be modified by choosing two points on the boundary of the rectangle.

\begin{theorem}
Let $R$ be a rectangle, with the origin $o$ as its center, and with side lengths $2$ and $2a$, where $a > 0$. Choose two points, $p_1$ and $p_2$ randomly and independently on the boundary of $R$, using uniform distribution. Let $R'$ denote the part of $R$, truncated by the segment $[p_1,p_2]$, containing $o$. For $S=3,4,5$, let $p_S(a)$ denote the probability that $R'$ has $S$ stable equilibrium points with respect to $o$. Then both $p_5(a)$, and $p_4(a)+p_5(a)$ are quasiconcave functions of $a$, with their unique maximum attained at $a=1$.
\end{theorem}

\begin{proof}
We present the proof only for $p_5(a)$, as the proof for $p_4(a)+p_5(a)$ is similar.

Let the vertices of $R$ be $(\pm 1, \pm a)$.
Since $p_5(a) = p_5 \left( \frac{1}{a} \right)$, we may assume that $0 < a \leq 1$.
Note that $R'$ has $5$ stable points if, and only if the midpoint of each side of $R$ belongs to $R'$, and the projection of $o$ onto the line of $[p_1,p_2]$ lies on $[p_1,p_2]$. An elementary consideration shows that the probability that each midpoint of $R$ belongs to $R'$ is $\frac{a}{2(1+a)^2}$. We compute the probability that $[p_1,p_2]$ contains a new stable point under the condition that each midpoint of $R$ belongs to $R'$.

Without loss of generality, we may assume that $p_1=(x,a)$, where $0 \leq x \leq 1$, and that $p_2 = (1,y)$, where $0 \leq y \leq a$.
Then there is an equilibrium point on $[p_1, p_2]$, if both angles $(p_1,p_2,o) \angle$ and $(p_2,p_1,o)$ are acute, or equivalently, if
$\langle -p_1, p_2-p_1 \rangle > 0$ and $\langle -p_2, p_1-p_2 \rangle > 0$.
This conditions can be written as $-x+x^2 -ax-a^2 > 0$ and $x^2 -ax+1-x > 0$.
From this, after some case analysis, we obtain that the conditions are equivalent to the following:
\begin{itemize}
\item if $0 < a \leq \frac{1}{2}$, then $0 < x < \frac{1-\sqrt{1-4a^2}}{2}$ and $0 \leq y \leq \frac{a^2+x-x^2}{a}$; or $\frac{1+\sqrt{1-4a^2}}{2} < x < 1-\frac{a^2}{4}$ and $0 \leq y \leq \frac{a^2+x-x^2}{a}$; or $1-\frac{a^2}{4} < x < 1$ and $0 \leq y \leq \frac{a-\sqrt{a^2-4+4x}}{2}$ or $\frac{a+\sqrt{a^2-4+4x}}{2} \leq y \leq \frac{a^2+x-x^2}{a}$.
\item if $\frac{1}{2} < a \leq 1$, then $0 < x < 1-\frac{a^2}{4}$ and $0 \leq y \leq \frac{a^2+x-x^2}{a}$; or $1-\frac{a^2}{4} < x \leq 1$ and
$0 \leq y \leq \frac{a-\sqrt{a^2-4+4x}}{2}$ or $\frac{a+\sqrt{a^2-4+4x}}{2} \leq y \leq \frac{a^2+x-x^2}{a}$.
\end{itemize}

Now $p_5(a)$ can be computed by simple integrations:
\[
p_5(a) =\left\{
\begin{array}{l}
\frac{-1+6a^2-a^4+(1-4a^2)^{\frac{3}{2}}}{12a(a+1)^2},\quad \hbox{if } 0 < a \leq \frac{1}{2};\\
\frac{-1+6a^2-a^4}{12a(a+1)^2},\quad \hbox{if } \frac{1}{2} < a \leq 1.
\end{array}
\right.
\]
From this, an elementary computation yields the assertion.
\end{proof}

\begin{figure}[ht]
\includegraphics[width=0.5\textwidth]{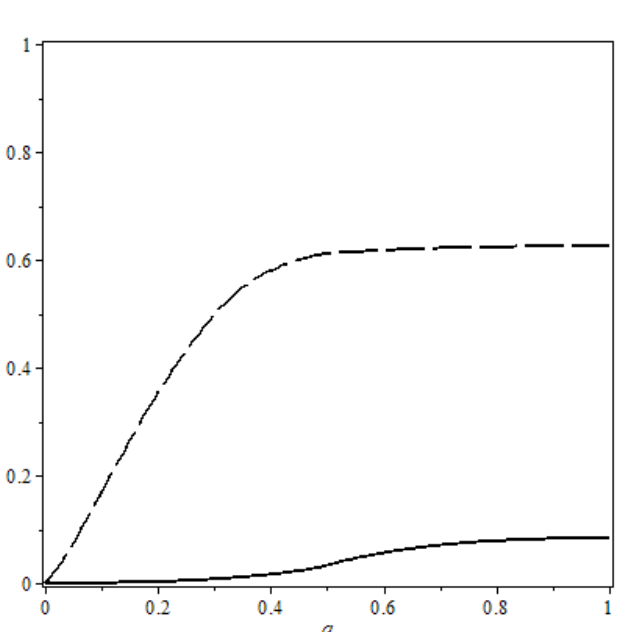}
\caption[]{Continuous line: the probability $p_5(a)$. Dashed line: the probability $p_4(a) + p_5(a)$. Both functions attain their maxima at $a=1$.}
\label{fig:rectangle_boundary}
\end{figure}

Figure~\ref{fig:rectangle_boundary} shows the probabilities $p_5(a)$ and $p_4(a)+p_5(a)$ as functions of $a$.

\section{Summary and Applications} \label{sec:sum}

In this paper we started to develop the theory of distance-driven flows as mathematical models of abrasion. The main thrust of the paper is
to identify geophysical shape descriptors which evolve either monotonically or in a quasiconcave or quasiconvex manner under some distance-driven flow.
Following the structure of the introduction, next we summarize our results by the geological shape descriptors.

\subsection{Summary of Results Grouped by Geological Shape Descriptors}
\begin{itemize}

\item{Axis ratios.}
We investigated axis ratios for ellipses and ellipsoids as initial conditions and showed that under radial flows given by a convex/concave function $f$ in Eq. (\ref{dist}),
axis ratios evolve monotonically but may achieve any limit as time approaches infinity. In case of orthogonal affinity (as a simple example
of a parallel flow) we showed that the axis ratio of a planar ellipse evolves in a quasiconcave manner, however, the larger axis ratio of
a tri-axial ellipsoid may have several temporal extrema.

\item{Isoperimetric ratio.}
We showed that in case of radial flows the evolution of the isoperimetric ratio is more complicated than that of the axis ratio (the former may exhibit more extrema
then the latter). On the other hand, in case of orthogonal affinity, the isoperimetric ratio always evolves in a quasiconcave manner so
it does not display more extrema than the evolution of the axis ratios.

\item{Number of static equilibrium points.}
We showed that both the number $S$ of stable and the number $U$ of unstable equilibrium points as well as the Morse-Smale complex ${\mathcal{ M}}(K)$ associated with
the gradient field are invariant under distance-driven flows.  According to Domokos et al. (2016a, 2016b), $\{S,U\}$ is called the primary equilibrium class of $K$ while ${\mathcal {M}}(K)$ uniquely defines the secondary and tertiary equilibrium classes of $K$.
Distance driven flows of the type in Eq. (\ref{dist}) can be interpreted as a continuous group acting on $\Re^n$ and based on
Theorem \ref{thm:constantT} and Corollary \ref{cor:constantT} the Morse-Smale complex is an invariant of these groups.
Since distance-driven flows have been introduced by Aristotle into the geometric theory of abrasion, based on the current
results we may call the primary, secondary and tertiary equilibrium classes the Aristotelian invariants of $K$.

\end{itemize}

\subsection{Questions and Conjectures}
One important conclusion from our results is that under one-parameter orthogonal affinity  $I(t)$ and $T(t)$ reach their respective minima simultaneously,
as $t$ approaches either zero or infinity. We also showed that $I(t)$ has a single maximum. While the maximum of $T(t)$ often
does not coincide with the maximum of $I(t)$, and $T(t)$ may have even several local maxima,  we still believe that the global trend of the two functions is related.
The one-parameter orthogonal affinity associates with each direction $\vv$ and each parameter value $t$ a real number $I_v(t) \in[0,1]$ and an integer $T_v(t)$.
One might try to formalize this relation by statistical methods.

Consider a convex polygon $P$ which is the convex hull of $m$ points chosen in a unit disk, independently and using uniform distribution.
Let $P(t)$ denote the image of $P$ under the orthogonal affinity defined by $(x,y) \mapsto (x,ty)$.
Let $T_m(t)$ be the expected value of the static equilibrium points of $P(t)$, over the family of all convex polygons with at most $n$ vertices,
using the probability distribution defined by the choice of $P$. Similarly, let $I_m(t)$ denote the expected value of the isoperimetric values of $P(t)$ using the same distribution.

\begin{conjecture}
Both $T_m(t)$ and $I_m(t)$ are quasiconcave functions for every $m \geq 3$, which attain their maxima at the same value of $t$.
\end{conjecture}

\subsection{Applications}

The results derived in this paper are of fundamental importance to explain and to interpret field  and laboratory data
if both collisional and frictional abrasion are significant. While the main focus of
our paper is theoretical, here we mention some immediate applications. In Miller et al. (2014) the shape evolution of pebbles was monitored
in the Bisley-Mameyes river system. In the field campaign several shape descriptors have been measured. Since the evolution of
pebbles was monitored from the original, fragmented shapes, axis ratios proved to be less reliable, however, the isoperimetric ratio was also measured.  One of the key observations of the paper is that shape evolution is caused partially
by collisions which dominate the initial phase of shape evolution, partially by friction in the second phase (Miller et al. 2014).
Similarly, in another field study (Szab\'o et al. 2013) along the Williams river, Australia, the combined effect of collisions
and friction has been pointed out.  So far, only the collisional part could be compared to mathematical
models, our current paper opens the possibility to study the combined action. In particular, if we accept orthogonal
affinity as a simple friction model, then Theorems \ref{thm:quasiconvex}, \ref{thm:basis}, \ref{thm:large} and \ref{thm:stochastic}, together with results on collisional abrasion (Bloore 1977; Firey 1974) lead to the following qualitative conclusions:

\begin{figure}[ht]
\includegraphics[width=0.9\textwidth]{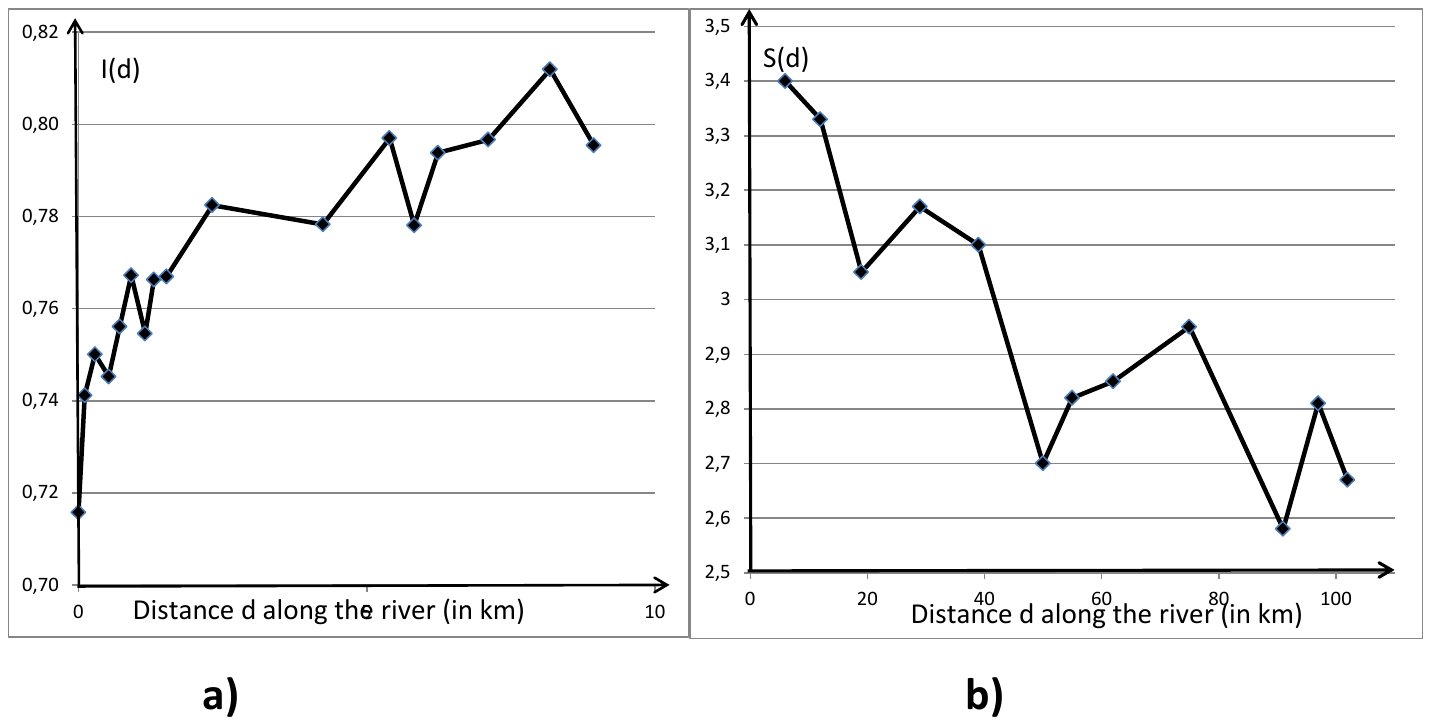}
\caption[]{a) Field data from Miller et al. )2014): Evolution of the isoperimetric ratio  in the Bisley-Mameyes river system , Puerto Rico. Observe initial sharp increase and subsequent saturation at $I=0.8$,
significantly below the maximum of $I=1$. b) Field data from Szab\'o et al. )2013):
Evolution of the number $S$ of stable static balance points. Observe decreasing trend with random fluctuations approaching the minimal value at $S=2$.}
\label{fig:data}
\end{figure}

\begin{itemize}
\item The isoperimetric ratio $I$ \emph{increases} under collisions but \emph{decreases} under friction.
We also note that under purely collisional abrasion the isoperimetric ratio $I$ increases monotonically
and saturates close to its maximum at $I=1$. If collisions
dominate the first phase and friction enters into the second phase then we expect $I$ to increase
initially sharply and subsequently to saturate/oscillate at a value significantly below the maximum of $I=1$.
\item The number of static balance points can be modeled by a random variable the expected value of which decreases \emph{both} under collision and friction, so in the
field data we expect a monotonically decreasing trend with random fluctuations, with either the stable or the unstable points approaching their minimal value  at $S=2$ or at $U=2$.
\end{itemize}
While the above conclusions are only qualitative, they are the first step towards the mathematical understanding of such diagrams. Figure \ref{fig:data} illustrates
that the theoretical predictions show a remarkably good match with the field data.

\section{Acknowledgements}
The authors gratefully acknowledge the support of the J\'anos Bolyai Research Scholarship of the Hungarian Academy of Sciences, and support from OTKA grant 119245. Furthermore, they express their gratitude to an unknown referee for making Theorems~\ref{thm:quasiconvex} and \ref{thm:basis} more general.

\end{document}